\documentclass[a4paper,11pt]{amsart}
\usepackage{amsfonts,amsmath,amsthm,amssymb}
\usepackage{enumitem}
\usepackage{mathrsfs}
\usepackage{booktabs}
\usepackage[margin=1in]{geometry}
\usepackage[colorlinks=true,linkcolor=blue,filecolor=blue,citecolor=red]{hyperref}
\usepackage{cleveref}

\theoremstyle{plain}
\newtheorem{theorem}{Theorem}[section]

\newtheorem{lemma}[theorem]{Lemma}

\newtheorem{proposition}[theorem]{Proposition}

\theoremstyle{definition}

\newtheorem{remark}[theorem]{Remark}

\begin{document}
\title[Analogue of Ramanujan's function $k(\tau)$ for the continued fraction $X(\tau)$]{Analogue of Ramanujan's function $k(\tau)$ for the continued fraction $X(\tau)$ of order six}
\author[Russelle Guadalupe, Victor Manuel Aricheta]{Russelle Guadalupe, Victor Manuel Aricheta}
\address{Institute of Mathematics, University of the Philippines-Diliman\\
Quezon City 1101, Philippines}
\email{rguadalupe@math.upd.edu.ph, vmaricheta@math.upd.edu.ph}

\renewcommand{\thefootnote}{}



\renewcommand{\thefootnote}{\arabic{footnote}}
\setcounter{footnote}{0}

\begin{abstract}
Motivated by the recent work of Park on the analogue of the Ramanujan's function $k(\tau)=r(\tau)r^2(2\tau)$ for the Ramanujan's cubic continued fraction, where $r(\tau)$ is the Rogers-Ramanujan continued fraction, we use the methods of Lee and Park to study the modularity and arithmetic of the function $w(\tau) = X(\tau)X(3\tau)$, which may be considered as an analogue of $k(\tau)$ for the continued fraction $X(\tau)$ of order six introduced by Vasuki, Bhaskar and Sharath. In particular, we show that $w(\tau)$ can be written in terms of the normalized generator $u(\tau)$ of the field of all modular functions on $\Gamma_0(18)$, and derive modular equations for $u(\tau)$ of smaller prime levels. We also express $j(d\tau)$ for $d\in\{1,2,3,6,9,18\}$ in terms of $u(\tau)$, where $j$ is the modular $j$-invariant.
\end{abstract}

\maketitle

\section{Introduction}\label{sec1}
For an element $\tau$ in the complex upper half-plane $\mathbb{H}$, define $q:= e^{2\pi i\tau}$. The Rogers-Ramanujan continued fraction $r(\tau)$ is defined by
\begin{align*}
r(\tau) = \dfrac{q^{1/5}}{1+\dfrac{q}{1+\dfrac{q^2}{1+\dfrac{q^3}{1+\cdots}}}} = q^{1/5}\prod_{n=1}^{\infty}\dfrac{(1-q^{5n-1})(1-q^{5n-4})}{(1-q^{5n-2})(1-q^{5n-3})}.
\end{align*}
Gee and Honsbeek \cite{geehons} computed the exact values of $r(\tau)$ at some imaginary quadratic points $\tau\in\mathbb{H}$ by proving that $r(\tau)$ generates the field of all modular functions on the congruence subgroup $\Gamma(5)$.
On the other hand, Ramanujan \cite[p. 326]{raman} showed that the function $k:=k(\tau)=r(\tau)r(2\tau)^2$ satisfies the relations
\begin{align}\label{eq11}
r^5(\tau) = k\left(\dfrac{1-k}{1+k}\right)^5, \qquad r^5(2\tau) = k^2\left(\dfrac{1+k}{1-k}\right).
\end{align}
Lee and Park \cite{leeparkk} proved that $k(\tau)$ generates the field of all modular functions on $\Gamma_1(10)$. Consequently, they showed that there is a modular equation for $k(\tau)$ of any level, and that $k(\tau)$ generates the ray class field modulo $10$ over an imaginary quadratic field $K$ for some suitable $\tau\in K\cap\mathbb{H}$. We refer the interested reader to \cite{coop, leeparkk, ye} for more details about the modular and arithmetic properties of $k(\tau)$. 

More recently, Park \cite{park} studied the modularity of the function $t(\tau) = C(\tau)C(2\tau)$, which is an analogue of $k(\tau)$ for Ramanujan's cubic continued fraction defined by
\begin{align*}
C(\tau) = \dfrac{q^{1/3}}{1+\dfrac{q+q^2}{1+\dfrac{q^2+q^4}{1+\dfrac{q^3+q^6}{1+\cdots}}}} = q^{1/3}\prod_{n=1}^{\infty}\dfrac{(1-q^{6n-1})(1-q^{6n-5})}{(1-q^{6n-3})^2}.
\end{align*}
Park established the following relations for $t:=t(\tau)$ given by 
\begin{align}\label{eq12}
C^3(\tau) = t(1-2t), \qquad C^3(2\tau) = \dfrac{t^2}{1-2t},
\end{align}
and proved that $t(\tau)$ generates the field of all modular functions on $\Gamma_0(12)$. Park also showed that there is a modular equation for $t(\tau)$ of any level, and that $t(\tau/2)$ generates the ray class field modulo $2$ over $K$ modulo $2$ for some suitable $\tau\in K\cap\mathbb{H}$.

Motivated by the work of Park, we explore in this paper the modularity and arithmetic of the function $w(\tau) := X(\tau)X(3\tau)$, which may be considered as an analogue of $k(\tau)$ for the continued fraction $X(\tau)$ of order six defined by (cf. \cite{park2})
\begin{align*}
X(\tau) &= \dfrac{q^{1/4}(1-q^2)}{1-q^{3/2}+\dfrac{(1-q^{1/2})(1-q^{7/2})}{q^{1/2}(1-q^{3/2})(1+q^3)+\dfrac{(1-q^{5/2})(1-q^{13/2})}{q^{3/2}(1-q^{3/2})(1-q^6)+\cdots}}}\\
&=q^{1/4}\prod_{n=1}^{\infty}\dfrac{(1-q^{6n-1})(1-q^{6n-5})}{(1-q^{6n-2})(1-q^{6n-4})}=\dfrac{\eta(\tau)\eta^2(6\tau)}{\eta^2(2\tau)\eta(3\tau)},
\end{align*}
where $\eta(\tau)=q^{1/24}\prod_{n=1}^\infty (1-q^n)$ is the Dedekind eta function. The function $X(\tau)$, which arises as a quotient of mock theta functions, is closely related to $C(\tau)$ via the well-known identity \cite[Thm. 6.9]{cooper}
\begin{align*}
\dfrac{1}{X^4(\tau)}-\dfrac{1}{C^3(\tau)}=1.
\end{align*}
Zagier \cite{zagier} found the following remarkable identities, labeled as Sequences \textbf{C} and \textbf{A}, that are valid in a neighborhood of $q=0$ given by (cf. \cite[p. 396, (6.57), (6.58)]{cooper})
\begin{align*}
\sum_{n=0}^\infty \sum_{\substack{j+k+\ell=n\\ j,k,\ell\geq 0}} \left(\dfrac{n!}{j!k!\ell!}\right)^2 X(\tau)^{4n} &=\dfrac{\eta^6(2\tau)\eta(3\tau)}{\eta^3(\tau)\eta^2(6\tau)}, \\
\sum_{n=0}^\infty \sum_{k=0}^n \dbinom{n}{k}^3 C(\tau)^{3n} &=\dfrac{\eta(2\tau)\eta^6(3\tau)}{\eta^2(\tau)\eta^3(6\tau)},
\end{align*}
which stem from the integral solutions of Beukers' differential equation and is connected with a proof of irrationality of $\zeta(3)$ by Ap\'{e}ry \cite{apery}. Additional identities for $X(\tau)$ were obtained by Vasuki, Bhaskar and Sharath \cite{vasuki}.

As a modular function, Lee and Park \cite{leeparkx} proved that $X^2(2\tau)$ generates the field of all modular functions on $\Gamma_0(12)$, and $X^2(\tau)$ generates the ray class field modulo $6$ over $K$ for some suitable $\tau\in K\cap\mathbb{H}$. Recently, Cho \cite{cho} proved their conjecture that $X(2\tau)$ in fact generates the field of all modular functions on
\begin{align*}
\left\langle \Gamma_1(12)\cap \Gamma^0(2), \begin{bmatrix}
-7 & -3\\ 12 & 5
\end{bmatrix}\right\rangle. 
\end{align*}

We begin our study on the modularity of $w(\tau)$ by showing that it can be written in terms of the $\eta$-quotient defined by 
\begin{align*}
u(\tau) := \dfrac{\eta(3\tau)\eta^3(18\tau)}{\eta(6\tau)\eta^3(9\tau)}, 
\end{align*}
which can be used to deduce our first two main results.

\begin{theorem}\label{thm11} We have the identities
\begin{align*}
X^4(\tau) = w(\tau)(1-3w(\tau)+3w(\tau)^2), \qquad X^4(3\tau) = \dfrac{w(\tau)^3}{1-3w(\tau)+3w(\tau)^2}.
\end{align*}
\end{theorem}

\begin{theorem}\label{thm12}
The functions $u(\tau)$ and $w(\tau)$ both generate the field of all modular functions on $\Gamma_0(18)$. 
\end{theorem} 

Theorem \ref{thm11} gives relations between $w(\tau), X(\tau)$ and $X(3\tau)$ similar to (\ref{eq11}) and (\ref{eq12}). In view of Theorem \ref{thm12}, $u(\tau)$ is a normalized generator for the field of all modular functions on $\Gamma_0(18)$ in the sense that the coefficient of $q^0$ in the $q$-expansion 
\begin{align*}
\dfrac{1}{u(\tau)} = \dfrac{1}{q}+q^2+q^5-q^8+O(q^9)
\end{align*}
of $1/u(\tau)$ is zero. In their recent paper \cite{hubsy}, Huber, Schultz and Ye obtained such generators for genus zero modular function fields, which are essential in constructing Ramanujan-Sato series for $1/\pi$. 

For each positive integer $n$, the modular equation for $u(\tau)$ of level $n$ is defined by $F_n(u(\tau), u(n\tau))=0$ for some irreducible polynomial $F_n(X,Y)\in\mathbb{C}[X,Y]$. The modular equation of any level can be explicitly computed, as the next result shows.

\begin{theorem}\label{thm13}
One can explicitly obtain modular equations for $u(\tau)$ and $w(\tau)$ of level $n$ for all positive integers $n$.
\end{theorem}

For example, the modular equations for $u(\tau)$ of levels two and three are, respectively:

\begin{align*}
u^2(\tau)-u(2\tau)+2u(\tau)u^2(2\tau)&=0,\\
u^3(\tau)-u(3\tau)+2u^3(\tau)u(3\tau)+u^2(3\tau)+4u^3(\tau)u^2(3\tau)-u^3(3\tau)&=0.
\end{align*}

Using these equations and the abovementioned relation between $u(\tau)$ and $w(\tau)$, we derive the modular equations for $w(\tau)$ of levels two and three (see Remark \ref{rem1}).

Finally, we show the following result that certain values of $w(\tau)$ generate ray class fields modulo $6$ over imaginary quadratic fields. 

\begin{theorem}\label{thm14}
Let $K$ be an imaginary quadratic field and $\tau\in K\cap\mathbb{H}$ be a root of a primitive equation $0=aX^2+bX+c$, with $a,b,c\in\mathbb{Z}$, such that the discriminant of $K$ is $b^2-4ac$. Suppose that $\gcd(a,6)=1$. Then $K(w(\tau/3))$ is the ray class field modulo $6$ over $K$.
\end{theorem}

We organize the paper as follows. We prove in Section \ref{sec2} Theorem \ref{thm11} by writing $w(\tau), X(\tau)$ and $X(3\tau)$ in terms of $u(\tau)$. In Section \ref{sec3}, we review some facts about modular functions, particularly $\eta$-quotients, on the congruence subgroup $\Gamma_0(N)$. In Section \ref{sec4}, we prove Theorem \ref{thm12}, and using an affine model for the modular curve $X_0(18)$, we deduce Theorem \ref{thm13}. We also establish modular equations for $u(\tau)$ of level $p$ for primes $p\leq 13$, and describe some properties of the modular equations for $u(\tau)$ of level $p$ for primes $p\geq 5$, including the Kronecker congruence relation. In Section \ref{sec5}, we describe certain arithmetic properties of $w(\tau)$ by establishing Theorem \ref{thm14}. We also show that $1/w(\tau)$ is an  algebraic integer for imaginary quadratic points $\tau\in\mathbb{H}$ by relating $u(\tau)$ with the modular $j$-invariant defined by 
\begin{align*}
j(\tau) = \eta^{-24}(\tau)\left(1+240\sum_{n=1}^\infty\dfrac{n^3q^n}{1-q^n}\right)^3.
\end{align*} 
Finally, in Section \ref{sec6}, we express $j(d\tau)$ in terms of $u(\tau)$ for $d\in\{2,3,6,9,18\}$. We use \textit{Mathematica} for our computations. 

\section{Proof of Theorem \ref{thm11}}\label{sec2}

\begin{proof}[Proof of Theorem \ref{thm11}]
By the theory of theta functions, we have the identities \cite[pp. 420-421]{cooper}
\begin{align*}
1+u = \dfrac{\eta_2^2\eta_3\eta_{18}}{\eta_1\eta_6\eta_9^2}, \qquad 1-2u=\dfrac{\eta_1^2\eta_{18}}{\eta_2\eta_9^2},
\end{align*}
and 
\begin{align*}
1+u^3 = \dfrac{\eta_6^5\eta_{18}}{\eta_3\eta_9^5}, \qquad 1-8u^3=\dfrac{\eta_3^8\eta_{18}^4}{\eta_6^4\eta_9^8},
\end{align*}
where $u = u(\tau)$ and $\eta_k := \eta(k\tau)$. Thus, if $e_1,e_2,e_3,e_6,e_9$ and $e_{18}$ are integers whose sum is zero, then there exist unique constants $a, b, c, d$ and $e$ such that
\begin{align*}
\eta_1^{e_1}\eta_2^{e_2}\eta_3^{e_3}\eta_6^{e_6}\eta_9^{e_9}\eta_{18}^{e_{18}}=u^a(1+u)^b(1-2u)^c(1+u^3)^d(1-8u^3)^e.
\end{align*} 
It follows from the above identities and the definitions of $X(\tau)$ and $w(\tau)$ that
\begin{align}
w(\tau) &= \dfrac{\eta(\tau)\eta^2(18\tau)}{\eta^2(2\tau)\eta(9\tau)} = \dfrac{u}{1+u},\label{eq21}\\
X^4(\tau) &= \dfrac{\eta^4(\tau)\eta^8(6\tau)}{\eta^8(2\tau)\eta^4(3\tau)} = \dfrac{u(1+u^3)}{(1+u)^4},\nonumber\\
X^4(3\tau) &= \dfrac{\eta^4(3\tau)\eta^8(18\tau)}{\eta^8(6\tau)\eta^4(9\tau)} = \dfrac{u^3}{1+u^3}.\nonumber
\end{align}
Hence, we arrive at
\begin{align*}
X^4(\tau) &= w(\tau)\cdot \dfrac{1+u^3}{(1+u)^3} = w(\tau)\cdot \left(1 - \dfrac{3u}{1+u}+\dfrac{3u^2}{(1+u)^2}\right) = w(\tau)(1-3w(\tau)+3w^2(\tau)),\\
X^4(3\tau) &= \dfrac{w^4(\tau)}{X^4(\tau)} = \dfrac{w^3(\tau)}{1-3w(\tau)+3w^2(\tau)},
\end{align*}
which completes the proof.
\end{proof}

\section{Modular functions on $\Gamma_0(N)$}\label{sec3}

We recall some facts about modular functions on the congruence subgroup
\begin{equation*}
\Gamma_0(N) := \left\lbrace\begin{bmatrix}
a & b\\ c& d
\end{bmatrix} \in \mbox{SL}_2(\mathbb{Z}) : c\equiv 0\pmod N\right\rbrace.
\end{equation*}
Any element of $\Gamma_0(N)$ acts on the extended upper half plane $\mathbb{H}^\ast := \mathbb{H}\cup \mathbb{Q}\cup \{\infty\}$ by a linear fractional transformation. We define the cusps of $\Gamma_0(N)$ to be the equivalence classes of $\mathbb{Q}\cup \{\infty\}$ under this action, and it is known that there are finitely many inequivalent cusps of $\Gamma_0(N)$. We define a modular function on $\Gamma_0(N)$ as a complex-valued function $f(\tau)$ on $\mathbb{H}$ such that $f(\tau)$ is meromorphic on $\mathbb{H}$, $f(\gamma\tau)=f(\tau)$ for all $\gamma\in \Gamma_0(N)$, and $f(\tau)$ is meromorphic at all cusps of $\Gamma_0(N)$. The last condition implies that for every cusp $r$ of $\Gamma_0(N)$ and an element $\gamma\in\mbox{SL}_2(\mathbb{Z})$ with $\gamma(r)=\infty$, $f\circ \gamma^{-1}$ has a $q$-expansion given by $(f\circ\gamma^{-1})(\tau) = \sum_{n\geq n_0} a_nq^{n/h}$ for some integers $h$ and $n_0$ with $a_{n_0}\neq 0$. We call $n_0$ the order of $f(\tau)$ at $r$, denoted by $\mbox{ord}_r f(\tau)$, and we say that $f(\tau)$ has a zero (resp., a pole) at $r$ if $\mbox{ord}_r f(\tau)$ is positive (resp., negative).

We denote by $A_0(\Gamma_0(N))$ the field of all modular functions on $\Gamma_0(N)$ and $A_0(\Gamma_0(N))_{\mathbb{Q}}$ its subfield consisting of all modular functions on $\Gamma_0(N)$ with rational $q$-expansion. We identify $A_0(\Gamma_0(N))$ with the field of all meromorphic functions on the modular curve $X_0(N):=\Gamma_0(N)\backslash \mathbb{H}^\ast$. If $f(\tau)\in A_0(\Gamma_0(N))$ has no zeros nor poles on $\mathbb{H}$, then the degree $[A_0(\Gamma_0(N)):\mathbb{C}(f(\tau))]$ is the total degree of poles of $f(\tau)$ given by $-\sum_r \mbox{ord}_r f(\tau)$, where the sum ranges over the inequivalent cusps $r$ of $\Gamma_0(N)$ for which $f(\tau)$ has a pole at $r$ (see \cite[Prop. 2.11]{shimura}). The following lemma gives a complete set of inequivalent cusps of $\Gamma_0(N)$.

\begin{lemma}[\cite{chokoop}]\label{lem31}
Let $a,a',c$ and $c'$ be integers with $\gcd(a,c)=\gcd(a',c')=1$. Let $S_{\Gamma_0(N)}$ be the set of all inequivalent cusps of $\Gamma_0(N)$. We denote $\pm1/0$ as $\infty$. Then 
\begin{enumerate}
\item $a/c$ and $a'/c'$ are equivalent over $\Gamma_0(N)$ if and only if there exist an integer $n$ and an element $\overline{s}\in (\mathbb{Z}/N\mathbb{Z})^\times$ such that $(a',c')\equiv(\overline{s}^{-1}a+nc,\overline{s}c)\pmod N$, and 
\item we have 
\begin{align*}
S_{\Gamma_0(N)} = \{a_{c,j}/c\in \mathbb{Q} : 0 < c \mid N, 0 < a_{c,j}\leq N, \gcd(a_{c,j},N)=1,\\ a_{c,j}= a_{c,j'}\stackrel{\text{def}}{\iff} a_{c,j}\equiv a_{c,j'}\pmod{\gcd(c,N/c)}\}.
\end{align*}
\end{enumerate}
\end{lemma}

An \textit{$\eta$-quotient} is a function of the form
\begin{equation*}
f(\tau) = \prod_{\delta\mid N} \eta(\delta\tau)^{r_{\delta}}
\end{equation*}
for some index set $\{r_\delta\in\mathbb{Z} : \delta\mid N\}$. The next two lemmas give necessary conditions on the modularity and behavior of an $\eta$-quotient at the cusps of $\Gamma_0(N)$.

\begin{lemma}\label{lem32}
Let $f(\tau) = \prod_{\delta\mid N} \eta(\delta\tau)^{r_{\delta}}$ be an $\eta$-quotient with $k= \frac{1}{2}\sum_{\delta\mid N} r_{\delta}\in\mathbb{Z}$ such that
\begin{equation*}
\sum_{\delta\mid N} \delta r_{\delta} \equiv 0\pmod{24}\quad\text{ and }\quad\sum_{\delta\mid N} \dfrac{N}{\delta}r_{\delta} \equiv 0\pmod {24}.
\end{equation*}
Then for all $[\begin{smallmatrix}
a & b\\ c & d
\end{smallmatrix}]\in \Gamma_0(N)$,
\begin{equation*}
f\left(\dfrac{a\tau+b}{c\tau+d}\right) = \left(\dfrac{(-1)^k\prod_{\delta\mid N} \delta^{r_\delta}}{d}\right)(c\tau+d)^kf(\tau).
\end{equation*}
\end{lemma} 

\begin{proof}
See \cite[Thm. 1.64]{ono} (or \cite{newman1, newman2}).
\end{proof}

\begin{lemma}\label{lem33}
Let $c, d$ and $N$ be positive integers with $d\mid N$ and $\gcd(c,d)=1$ and let $f(\tau) = \prod_{\delta\mid N} \eta(\delta\tau)^{r_{\delta}}$ be an $\eta$-quotient satisfying the conditions of Lemma \ref{lem32}. Then the order of vanishing of $f(\tau)$ at the cusp $c/d$ is 
\begin{align*}
\dfrac{N}{24d\gcd(d,\frac{N}{d})}\sum_{\delta\mid N} \gcd(d,\delta)^2\cdot\dfrac{r_{\delta}}{\delta}.
\end{align*}
\end{lemma}

\begin{proof}
See \cite[Thm. 1.65]{ono} (or \cite{ligo}).
\end{proof}

\section{Modularity of $u(\tau)$ and $w(\tau)$}\label{sec4}

In this section, we prove Theorem \ref{thm12} and deduce that the modular equations for $u(\tau)$ and $w(\tau)$ exist at any level.

\begin{proof}[Proof of Theorem \ref{thm12}]
By Lemma \ref{lem32}, $u(\tau)$ is a modular function on $\Gamma_0(18)$. We now take
\[S_{\Gamma_0(18)} = \{\infty, 0,1/2,1/3,2/3,1/6,5/6,1/9\}\] as the set of inequivalent cusps of $\Gamma_0(18)$ by Lemma \ref{lem31}. We compute the order of $u(\tau)$ at each element of $S_{\Gamma_0(18)}$ using Lemma \ref{lem33}, as shown in Table \ref{tbl41}.
\begin{table}[h]
\caption{The orders of $u(\tau)$ at the cusps of $\Gamma_0(18)$}\label{tbl41}
\begin{tabular}{@{}ccccccc@{}}
\hline
cusp $r$ & $\infty$  & $0$ & $1/2$  & $1/3, 2/3$ & $1/6, 5/6$ & $1/9$ \\
\hline
$\mbox{ord}_r u(\tau)$ & $1$ & $0$ & $0$ & $0$ & $0$ & $-1$\\
\hline
\end{tabular}
\end{table}

Thus, $u(\tau)$ has a simple pole at $1/9$ and a simple zero at $\infty$. Hence, $A_0(\Gamma_0(18)) = \mathbb{C}(u(\tau))$, and in view of (\ref{eq21}), we have $A_0(\Gamma_0(18)) = \mathbb{C}(w(\tau))$.
\end{proof}

As the modular curve $X_0(18)$ has genus zero, the next lemma shows that there is an affine plane model for $X_0(18n)$ over $\mathbb{Q}$ defined by $u(\tau)$. We omit the proof as it is similar to that of \cite[Lem. 3.4]{leeparkx}.

\begin{lemma}\label{lem41}
For all positive integers $n$, we have $A_0(\Gamma_0(18n))_{\mathbb{Q}}= \mathbb{Q}(u(\tau),u(n\tau))$.
\end{lemma}

We next observe the behavior of $u(\tau)$ at all cusps $r\in \mathbb{Q}\cup\{\infty\}$.

\begin{lemma}\label{lem42}
Let $a, c, a', c'$ and $n$ be integers with $n$ positive. Then:
\begin{enumerate}
\item The modular function $u(\tau)$ has a pole at $a/c\in \mathbb{Q}\cup\{\infty\}$ if and only if $\gcd(a,c) = 1$ and $c\equiv 9\pmod{18}$.
\item The modular function $u(n\tau)$ has a pole at $a'/c'\in \mathbb{Q}\cup\{\infty\}$ if and only if there exist integers $a$ and $c$ such that $a/c=na'/c', \gcd(a,c)=1$ and $c\equiv 9\pmod{18}$.
\item The modular function $u(\tau)$ has a zero at $a/c\in \mathbb{Q}\cup\{\infty\}$ if and only if $\gcd(a,c) = 1$ and $18\mid c$.
\item The modular function $u(n\tau)$ has a zero at $a'/c'\in \mathbb{Q}\cup\{\infty\}$ if and only if there exist integers $a$ and $c$ such that $a/c=na'/c', \gcd(a,c)=1$ and $18\mid c$.
\end{enumerate}
\end{lemma}

\begin{proof}
Using Table \ref{tbl41}, we know that $u(\tau)$ has a simple pole (resp., zero) at $a/c\in \mathbb{Q}\cup\{\infty\}$ which is equivalent to $1/9$ (resp., $\infty$) on $\Gamma_0(18)$. By Lemma \ref{lem31}, $a/c$ is equivalent to $1/9$ on $\Gamma_0(18)$ if and only if $(a,c)\equiv (s^{-1}+9n, 9s)\pmod{18}$ for some integer $n$ and $s\in(\mathbb{Z}/18\mathbb{Z})^\times$. This implies that $\gcd(a,c)=1$ and $c\equiv 9\pmod{18}$, proving (1). On the other hand, $a/c$ is equivalent to $\infty$ on $\Gamma_0(18)$ if and only if $(a,c)\equiv (s^{-1}, 0)\pmod{18}$ for some integer $n$ and $s\in(\mathbb{Z}/18\mathbb{Z})^\times$. This implies that $\gcd(a,c)=1$ and $18\mid c$, proving (3). Statements (2) and (4) immediately follow from (1) and (3).
\end{proof}

Ishida and Ishii \cite{ishdai} proved the following result, which can be used to describe the coefficients of the modular equation between two modular functions that are generators of the field $A_0(\Gamma')$, where $\Gamma'$ is a congruence subgroup of $\mbox{SL}_2(\mathbb{Z})$.

\begin{proposition}\label{prop43}
Let $\Gamma'$ be a congruence subgroup of $\mbox{SL}_2(\mathbb{Z})$, and $f_1(\tau)$ and $f_2(\tau)$ be nonconstant functions such that $A_0(\Gamma')=\mathbb{C}(f_1(\tau),f_2(\tau))$. For $k\in\{1,2\}$, let $d_k$ be the total degree of poles of $f_k(\tau)$. Let 
\begin{equation*}
F(X,Y)=\sum_{\substack{0\leq i\leq d_2\\0\leq j\leq d_1}}C_{i,j}X^iY^j\in \mathbb{C}[X,Y]
\end{equation*}
satisfy $F(f_1(\tau),f_2(\tau))=0$. Let $S_{\Gamma'}$ be the set of all the inequivalent cusps of $\Gamma'$, and for $k\in\{1,2\}$, define 
\begin{equation*}
S_{k,0}= \{r\in S_{\Gamma'} : f_k(\tau)\text{ has a zero at }r\}\text{ and }S_{k,\infty}= \{r\in S_{\Gamma'} : f_k(\tau)\text{ has a pole at }r\}.
\end{equation*}
Further let
\begin{equation*}
a = -\sum_{r\in S_{1,\infty}\cap S_{2,0}}\mbox{ord}_r f_1(\tau)\quad\text{ and }\quad b = \sum_{r\in S_{1,0}\cap S_{2,0}}\mbox{ord}_r f_1(\tau)
\end{equation*}
with the convention that $a=0$ (resp. $b=0$) whenever $S_{1,\infty}\cap S_{2,0}$ (resp. $S_{1,0}\cap S_{2,0}$) is empty. Then
\begin{enumerate}
\item $C_{d_2,a}\neq 0$ and if, in addition, $S_{1,\infty}\subseteq S_{2,\infty}\cup S_{2,0}$, then $C_{d_2,j}=0$ for all $j\neq a$;
\item $C_{0,b}\neq 0$ and if, in addition, $S_{1,0}\subseteq S_{2,\infty}\cup S_{2,0}$, then $C_{0,j}=0$ for all $j\neq b$.
\end{enumerate}
By interchanging the roles of $f_1(\tau)$ and $f_2(\tau)$, one may obtain properties analogous to (1) - (2).
\end{proposition} 

\begin{proof}
See \cite[Lem. 3 and Lem. 6]{ishdai} (or \cite[Thm. 3.3]{chokimkoo}).
\end{proof}

We can now apply Lemma \ref{lem41} and Proposition \ref{prop43} to prove Theorem \ref{thm13}. We remark that the modular equations for $u(\tau)$ and $w(\tau)$ serve as affine plane models for $X_0(18n)$ over $\mathbb{Q}$.

\begin{proof}[Proof of Theorem \ref{thm13}]
We follow the argument used to prove \cite[Thm. 3.2]{leeparkx}. Let $\mathbb{C}(f_1(\tau),f_2(\tau))$ be the field of all modular functions on some congruence subgroup, where $f_1(\tau)$ and $f_2(\tau)$ are nonconstant functions. For $j\in \{1,2\}$, the degree $[\mathbb{C}(f_1(\tau),f_2(\tau)):\mathbb{C}(f_j(\tau))]$ is equal to the total degree $d_j$ of poles of $f_j(\tau)$. Thus, by \cite[Lem. 3.1]{chokimkoo}, there is a unique polynomial $\Phi(X,Y)\in \mathbb{C}[X,Y]$ such that $\Phi(X,f_2(\tau))$ (resp. $\Phi(f_1(\tau),Y)$) is the minimal polynomial of $f_1(\tau)$ (resp. $f_2(\tau)$) over $\mathbb{C}(f_2(\tau))$ (resp. $\mathbb{C}(f_1(\tau))$) of degree $d_2$ (resp. $d_1$). As $u(\tau)$ has rational $q$-expansion, setting $f_1(\tau) = u(\tau)$ and $f_2(\tau) = u(n\tau)$ and using Proposition \ref{prop43}, we see that there is a polynomial $F_n(X,Y)\in \mathbb{Q}[X,Y]$ such that $F_n(u(\tau),u(n\tau))=0$ for all positive integers $n$ with $\deg_X F_n(X,Y) = d_2$ and $\deg_Y F_n(X,Y)=d_1$. The equation $F_n(u(\tau),u(n\tau))=0$ is the modular equation for $u(\tau)$ of level $n$ for all positive integers $n$. We infer from (\ref{eq21}) that for all positive integers $n$,
\begin{align*}
F_n\left(\dfrac{w(\tau)}{1-w(\tau)},\dfrac{w(n\tau)}{1-w(n\tau)}\right)=0.
\end{align*}
Consider the polynomial
\begin{align*}
\mathcal{G}_n(X,Y) =(1-X)^{d_2}(1-Y)^{d_1}F_n\left(\dfrac{X}{1-X},\dfrac{Y}{1-Y}\right),
\end{align*}
so that $\mathcal{G}_n(w(\tau),w(n\tau))=0$. We choose an irreducible factor $G_n(X,Y)\in\mathbb{Q}[X,Y]$ of $\mathcal{G}_n(X,Y)$ with $G_n(w(\tau),w(n\tau))=0$ as follows: after plugging in the $q$-expansions of $w(\tau)$ and $w(n\tau)$ into the irreducible factors of $\mathcal{G}_n(X,Y)$, we find which among them vanishes in a neighborhood of $q=0$. Hence, we get that $G_n(w(\tau),w(n\tau))=0$ is the modular equation for $w(\tau)$ of level $n$ for all positive integers $n$. 
\end{proof}

We apply Theorem \ref{thm13} to find explicitly the modular equations for $u(\tau)$ of levels two and three, and obtain information about the coefficients of the modular equations for $u(\tau)$ of level $p$ for primes $p\geq 5$.

\begin{theorem}[Modular equation of level two]\label{thm44} We have 
\begin{equation*}
u^2(\tau)-u(2\tau)+2u(\tau)u^2(2\tau)=0.
\end{equation*}
\end{theorem}

\begin{proof}
In view of Lemmata \ref{lem31}, \ref{lem41} and \ref{lem42}, we work on the relevant cusps of $\Gamma_0(36)$:  $\infty, 1/9$ and $1/18$. By Lemma \ref{lem33}, we know that $u(\tau)$ has a double pole at $1/9$ and two simple zeros at $\infty$ and $1/18$. Also, $u(2\tau)$ has a double zero at $\infty$ and two simple poles at $1/9$ and $1/18$. Thus, the total degrees of poles of both $u(\tau)$ and $u(2\tau)$ are $2$, so by Proposition \ref{prop43}, there is a polynomial 
\begin{align*}
F_2(X,Y) = \sum_{0\leq i, j\leq 2}C_{i,j}X^iY^j \in\mathbb{C}[X,Y]
\end{align*}
such that $F_2(u(\tau),u(2\tau))=0$. Moreover, we get $C_{2,0}\neq 0, C_{0,1}\neq 0$ and $C_{2,1}=C_{2,2}=C_{0,0}=C_{0,2}=0$. Switching the roles of $u(\tau)$ and $u(2\tau)$ gives $C_{1,2}\neq 0$ and $C_{1,0}=0$. Using the $q$-expansion of $u(\tau)$ given by
\begin{align}
u(\tau) = q-q^4+2q^{10}-2q^{13}-q^{16}+4q^{19}+O(q^{20}),\label{eq41}
\end{align}
we can take $C_{0,1}=-1$. Hence, we obtain $F_2(X,Y) = X^2-Y+2XY^2$ and the desired identity follows.
\end{proof}

\begin{theorem}[Modular equation of level three]\label{thm45} We have 
\begin{equation*}
u^3(\tau)-u(3\tau)+2u^3(\tau)u(3\tau)+u^2(3\tau)+4u^3(\tau)u^2(3\tau)-u^3(3\tau)=0.
\end{equation*}
\end{theorem}

\begin{proof}
Using Lemmata \ref{lem31}, \ref{lem41} and \ref{lem42}, we work on the following cusps of $\Gamma_0(54)$:  $\infty, 1/9,5/9,\\1/18,5/18$ and $1/27$. We infer from Lemma \ref{lem33} that $u(\tau)$ has three simple poles at $1/9, 5/9$ and $1/27$, and three simple zeros at $1/18, 5/18$ and $\infty$. Also, $u(3\tau)$ has a triple pole at $1/27$ and a triple zero at $\infty$. We see that the total degrees of poles of both $u(\tau)$ and $u(3\tau)$ are $3$, so Proposition \ref{prop43} says that there is a polynomial 
\begin{align*}
F_3(X,Y) = \sum_{0\leq i, j\leq 3}C_{i,j}X^iY^j \in\mathbb{C}[X,Y]
\end{align*}
such that $F_3(u(\tau),u(3\tau))=0$. We have $C_{3,0}\neq 0$ and $C_{0,1}\neq 0$. Interchanging the roles of $u(\tau)$ and $u(3\tau)$, we deduce that $C_{0,3}\neq 0$ and $C_{i,3}=C_{j,0}=0$ for $i\in\{1,2,3\}$ and $j\in \{0,1,2\}$. Using the $q$-expansion (\ref{eq41}), we can take $C_{0,1}=-1$. Hence, we arrive at $F_3(X,Y)=X^3-Y+2X^3Y+Y^2+4X^3Y^2-Y^3$ and the desired identity follows.
\end{proof}

\begin{theorem}\label{thm46}
Let $p\geq 5$ be a prime and $0=F_p(X,Y) = \sum_{0\leq i, j\leq p+1} C_{i,j} X^iY^j\in \mathbb{Q}[X,Y]$ be the modular equation of $u(\tau)$ of level $p$. Then
\begin{enumerate}
\item $C_{p+1,0}\neq 0$ and $C_{0,p+1}\neq 0$,
\item $C_{p+1,j}  = C_{j,p+1} = 0$ for $j\in \{1,2,\ldots, p+1\}$, and
\item $C_{0,j}  = C_{j,0} = 0$ for $j\in \{0,1,\ldots, p\}$.
\end{enumerate}
\end{theorem}

\begin{proof}
Using Lemmata \ref{lem31}, \ref{lem41} and \ref{lem42}, we work on the following cusps of $\Gamma_0(18p)$:  $\infty, 1/9,1/18$ and $1/9p$. We know from Lemma \ref{lem33} that $\mbox{ord}_{1/9} u(\tau)=\mbox{ord}_{1/9p} u(p\tau)=-p$ and $\mbox{ord}_{1/9p} u(\tau)=\mbox{ord}_{1/9} u(p\tau)=-1$. Thus, the total degrees of poles of both $u(\tau)$ and $u(p\tau)$ are $p+1$, so by Proposition \ref{prop43}, there is a polynomial
\begin{align*}
F_p(X,Y) = \sum_{0\leq i, j\leq p+1}C_{i,j}X^iY^j \in\mathbb{Q}[X,Y]
\end{align*}
such that $F_p(u(\tau),u(p\tau))=0$. We also have $C_{p+1,0}\neq 0$ and $C_{p+1,j}=0$ for $j\in \{1,2,\ldots,p+1\}$. As $\mbox{ord}_{1/18} u(\tau)=p$ and $\mbox{ord}_{\infty} u(\tau)=1$ by Lemma \ref{lem33}, we deduce that $C_{0,p+1}\neq 0$ and $C_{0,j}=0$ for $j\in \{0,1,\ldots,p\}$. Switching the roles of $u(\tau)$ and $u(p\tau)$, we obtain $C_{j,p+1}=0$ for $j\in \{1,2,\ldots,p+1\}$ and $C_{j,0}=0$ for $j\in \{0,1,\ldots,p\}$.
\end{proof}

The modular equations $F_p(X,Y)=0$ for $u(\tau)$ of prime level $p\leq 13$ are shown in Table \ref{tbl42}. We observe that for these primes $p$, the polynomial $F_p(X,Y)$ satisfy the \textit{Kronecker congruence relation}
\begin{equation*}
F_p(X,Y)\equiv (X^p-Y)(X-Y^p)\pmod{p}.
\end{equation*}

\begin{table}[h]
\caption{Modular equations $F_p(X,Y)=0$ for the function $u(\tau)$ of level $p$ for primes $p\leq 13$}\label{tbl42}%
{\small \begin{tabular}{@{}ll@{}}
\toprule
$p$ & $F_p(X,Y)$ \\
\midrule
$2$ & $X^2-Y+2XY^2$ \\
\midrule
$3$ & $X^3 - Y + 2 X^3 Y + Y^2 + 4 X^3 Y^2 - Y^3$\\
\midrule
$5$ & $(X^5 - Y) (X - Y^5) - 5 X Y (-X^3 - X Y + 2 X^4 Y + 4 X^2 Y^2 - Y^3 - 4 X^3 Y^3 + 
2 X Y^4 + 3 X^4 Y^4)$\\
\midrule
$7$ & $\begin{aligned}
&(X^7 - Y) (X - Y^7) - 7 X Y (-X^3 - 4 X^5 Y + 8 X^4 Y^2 - Y^3 - 3 X^3 Y^3 + 8 X^6 Y^3 + 8 X^2 Y^4 \\
&- 4 X Y^5 + 8 X^3 Y^6 + 9 X^6 Y^6)	
\end{aligned}$\\
\midrule
$11$ & $\begin{aligned}
&(X^{11} - Y) (X - Y^{11}) - 11 X Y (-X^3 + 2 X^6 + 2 X^9 + X Y - 6 X^4 Y - 17 X^7 Y - 4 X^{10} Y\\
&+2 X^2 Y^2 - 14 X^5 Y^2 + 4 X^8 Y^2 - Y^3 + 8 X^3 Y^3 -14 X^6 Y^3 - 68 X^9 Y^3 - 6 X Y^4 + 62 X^4 Y^4\\
&+ 28 X^7 Y^4 +	32 X^{10} Y^4 - 14 X^2 Y^5 + 14 X^5 Y^5 + 112 X^8 Y^5 + 2 Y^6 -14 X^3 Y^6 + 248 X^6 Y^6\\
&+ 192 X^9 Y^6 - 17 X Y^7 + 28 X^4 Y^7 + 128 X^7 Y^7 + 128 X^{10} Y^7 + 4 X^2 Y^8 + 112 X^5 Y^8 +128 X^8 Y^8 \\
&+ 2 Y^9 - 68 X^3 Y^9 + 192 X^6 Y^9 + 256 X^9 Y^9 - 4 X Y^{10} + 32 X^4 Y^{10} + 128 X^7 Y^{10} + 93 X^{10} Y^{10})
\end{aligned}$\\
\midrule
$13$ & $\begin{aligned}
&(X^{13} - Y) (X - Y^{13}) - 13 X Y (-X^3 + 3 X^6 + X^9 - 2 X^{12} + 2 X^2 Y - 12 X^5 Y - 30 X^8 Y\\
&- 23 X^{11} Y + 2 X Y^2 - 3 X^4 Y^2 - 93 X^7 Y^2 - 60 X^{10} Y^2 - Y^3 + 13 X^3 Y^3 - 33 X^6 Y^3 - 132 X^9 Y^3\\
&- 8 X^{12} Y^3 - 3 X^2 Y^4 - 9 X^5 Y^4 + 388 X^8 Y^4 + 240 X^{11} Y^4 - 12 X Y^5 - 9 X^4 Y^5 + 516 X^7 Y^5\\
&+ 744 X^{10} Y^5 + 3 Y^6 - 33 X^3 Y^6 + 300 X^6 Y^6 + 264 X^9 Y^6 + 192 X^{12} Y^6 - 93 X^2 Y^7 + 516 X^5 Y^7\\
&+ 72 X^8 Y^7 - 768 X^{11} Y^7 - 30 X Y^8 + 388 X^4 Y^8 + 72 X^7 Y^8 - 192 X^{10} Y^8 + Y^9 - 132 X^3 Y^9\\
&+ 264 X^6 Y^9 + 832 X^9 Y^9 + 512 X^{12} Y^9 - 60 X^2 Y^{10} + 744 X^5 Y^{10} - 192 X^8 Y^{10} - 1024 X^{11} Y^{10}\\
&- 23 X Y^{11} + 240 X^4 Y^{11} - 768 X^7 Y^{11} - 1024 X^{10} Y^{11} - 2 Y^{12} -8 X^3 Y^{12} + 192 X^6 Y^{12} \\
&+ 512 X^9 Y^{12} + 315 X^{12} Y^{12})
\end{aligned}$\\
\bottomrule
\end{tabular}}
\end{table}

\begin{remark}\label{rem1}
Following the proof of Theorem \ref{thm13}, one can find the corresponding modular equations for $w(\tau)$ of prime level $p\leq 13$ using $F_p(X,Y)=0$. For instance, we deduce the modular equations for $w(\tau)$ of levels two and three:
\begin{align*}
&w(\tau)^2-w(2\tau)+2w(\tau)w(2\tau)-3w(\tau)^2w(2\tau)+w(2\tau)^2=0,\\
w(\tau)^3-w(3\tau)&+3w(\tau)w(3\tau)-3w(\tau)^2w(3\tau)+3w(3\tau)^2-9w(\tau)w(3\tau)^2\\
&+9w(\tau)^2w(3\tau)^2-3w(3\tau)^3+9w(\tau)w(3\tau)^3-9w(\tau)^2w(3\tau)^3=0.
\end{align*}
\end{remark}

We are in a position to describe the properties of the modular equation for $u(\tau)$ of levels coprime to $6$. For an integer $a$ coprime to $6$, we define matrices $\sigma_a\in \mbox{SL}_2(\mathbb{Z})$ such that $\sigma_a\equiv [\begin{smallmatrix}
a^{-1} & 0\\ 0 & a
\end{smallmatrix}]\pmod{18}$. Then $\sigma_a\in \Gamma_0(18)$, and we may take
\begin{equation*}
\sigma_{\pm 1} =\pm\begin{bmatrix}
1 & 0\\ 0 & 1
\end{bmatrix},\qquad\sigma_{\pm 5}=\pm\begin{bmatrix}
-7 & -72\\18 & 185 
\end{bmatrix}\quad\text{ and }\quad\sigma_{\pm 7}=\pm\begin{bmatrix}
-5 & -72\\18 & 259
\end{bmatrix}.
\end{equation*}

For any integer $n$ coprime to $6$, we have the disjoint union 
\begin{equation*}
\Gamma_0(18)\begin{bmatrix}
1 & 0\\0 & n
\end{bmatrix}\Gamma_0(18) = \bigsqcup_{0<a\mid n}\bigsqcup_{\substack{0\leq b < n/a\\\gcd(a,b,n/a)=1}}\Gamma_0(18)\sigma_a\begin{bmatrix}
a & b\\0 & \frac{n}{a}
\end{bmatrix}
\end{equation*}
and $\Gamma_0(18)\backslash \Gamma_0(18)[\begin{smallmatrix}
1 & 0\\0 & n
\end{smallmatrix}]\Gamma_0(18)$ has $\psi(n):=n\prod_{p\mid n}(1+1/p)$ elements (see \cite[Prop. 3.36]{shimura}). We now define the polynomial
\begin{equation*}
\Phi_n(X,\tau):= \prod_{0<a\mid n}\prod_{\substack{0\leq b < n/a\\\gcd(a,b,n/a)=1}}(X-(u\circ\alpha_{a,b})(\tau))
\end{equation*}
where $\alpha_{a,b}:=\sigma_a[\begin{smallmatrix}
a & b\\0 & n/a
\end{smallmatrix}]$. Note that the coefficients of $\Phi_n(X,\tau)$ are elementary symmetric functions of $u\circ\alpha_{a,b}$, so these are invariant under the action of $\Gamma_0(18)$ and thus are in $A_0(\Gamma_0(18))=\mathbb{C}(u(\tau))$. We may then write $\Phi_n(X,\tau)$ as $\Phi_n(X,u(\tau))\in \mathbb{C}(u(\tau))[X]$. Since $\alpha_{1,0}=\sigma_1[\begin{smallmatrix}
1 & 0\\0 & n
\end{smallmatrix}]$ and $(u\circ\alpha_{1,0})(\tau) = u(\tau/n)$ is a zero of $\Phi_n(X,u(\tau))$, we have $\Phi_n(u(\tau/n),u(\tau))=0$, which can be seen as the modular equation for $u(\tau)$ of level $n$. The following theorem shows several properties of the polynomial $\Phi_n(X,Y)$ when $\gcd(n,6)=1$. We omit the proof as it is similar to that of \cite[Thm. 3.9]{leeparkx}.

\begin{theorem}\label{thm47}
Let $n$ be a positive integer coprime to $6$ and let $\Phi_n(X,Y)$ be the polynomial defined as above. Then
\begin{enumerate}
\item $\Phi_n(X,Y)\in \mathbb{Z}[X,Y]$ and $\deg_X\Phi_n(X,Y)=\deg_Y\Phi_n(X,Y)=\psi(n)$.
\item $\Phi_n(X,Y)$ is irreducible both as a polynomial in $X$ over $\mathbb{C}(Y)$ and as a polynomial in $Y$ over $\mathbb{C}(X)$.
\item $\Phi_n(X,Y)=\Phi_n(Y,X)$.
\item (Kronecker congruence relation) If $p\geq 5$ is a prime, then 
\begin{equation*}
\Phi_p(X,Y)\equiv (X^p-Y)(X-Y^p)\pmod{p\mathbb{Z}[X,Y]}.
\end{equation*} 
\end{enumerate}
\end{theorem}

\section{Arithmetic of $w(\tau)$}\label{sec5}

In this section, we provide some arithmetic properties of $w(\tau)$. We start with the following lemma, which asserts that ray class fields over an imaginary quadratic field can be generated by singular values of a suitable modular function on some congruence subgroup between $\Gamma(N)$ and $\Gamma_1(N)$. For $N\geq 1$, we denote $\Gamma(N)$ and $\Gamma_1(N)$ the set of all matrices $[\begin{smallmatrix}
a & b\\c & d
\end{smallmatrix}]\in\mbox{SL}_2(\mathbb{Z})$ congruent to $[\begin{smallmatrix}
1 & 0\\0 & 1
\end{smallmatrix}]$ and $[\begin{smallmatrix}
1 & \ast\\0 & 1
\end{smallmatrix}]$, respectively, modulo $N$.

\begin{lemma}\label{lem51}
Let $K$ be an imaginary quadratic field with discriminant $d_K$ and $\tau\in K\cap\mathbb{H}$ be a root of a primitive equation $aX^2+bX+c=0$ such that $b^2-4ac=d_K$ and $a,b,c\in\mathbb{Z}$. Let $\Gamma'$ be a congruence subgroup such that $\Gamma(N)\subset \Gamma'\subset \Gamma_1(N)$. Suppose that $\gcd(a,N)=1$. Then the field generated over $K$ by all values of $h(\tau)$, where $h\in A_0(\Gamma')_\mathbb{Q}$ is defined and finite at $\tau$, is the ray class field modulo $N$ over $K$.
\end{lemma}

\begin{proof}
See \cite[Cor. 5.2]{chokoo}.	
\end{proof}

We now use Lemma \ref{lem51} to prove Theorem \ref{thm14}.

\begin{proof}[Proof of Theorem \ref{thm14}] Since $w(\tau)$ has rational $q$-expansion given by 
\begin{align*}
	w(\tau) = q-q^2+q^3-2q^4+3q^5-4q^6+5q^7+O(q^8),
\end{align*}
we have $A_0(\Gamma_0(18))_{\mathbb{Q}}=\mathbb{Q}(w(\tau))$ in view of Theorem \ref{thm12}. For a subgroup $\Gamma'$ of $\mbox{SL}_2(\mathbb{Z})$, define $\overline{\Gamma'} := \Gamma'/\langle -I\rangle$. Since 
\begin{equation*}
	\begin{bmatrix}
		1 & 0\\ 0 & 3
	\end{bmatrix}^{-1}\overline{\Gamma_0(18)}\begin{bmatrix}
		1 & 0\\ 0 & 3
	\end{bmatrix} = \overline{\Gamma_0(6)}\cap \overline{\Gamma(3)}=\overline{\Gamma_1(6)}\cap \overline{\Gamma(3)},
\end{equation*}
we see that $\mathbb{Q}(w(\tau/3)) = A_0(\Gamma_1(6)\cap \Gamma(3))_{\mathbb{Q}}$. Let $K$ be an imaginary quadratic field and $\tau\in K\cap\mathbb{H}$ be a root of a primitive equation $0=aX^2+bX+c\in\mathbb{Z}[X]$ such that the discriminant of $K$ is $b^2-4ac$. Since $\gcd(6,a)=1$, $\Gamma(6)\subset \Gamma_1(6)\cap \Gamma(3)\subset\Gamma_1(6)$ and $w(\tau/3)$ is well-defined and finite, by Lemma \ref{lem51} with $N=6$, we see that the ray class field modulo $6$ over $K$ is $K(w(\tau/3))$.
\end{proof}

We next study the singular values of $1/w(\tau)$ at imaginary quadratic points $\tau\in\mathbb{H}$ by relating $u(\tau)$ with the modular $j$-invariant $j(\tau)$ mentioned in the introduction.

\begin{theorem}\label{thm52}
Let $K$ be an imaginary quadratic field and $\tau\in K\cap\mathbb{H}$. Then $1/w(\tau)$ is an algebraic integer.
\end{theorem}

\begin{proof}
We first remark that the $\eta$-quotient $f(\tau) := \eta^{24}(\tau)\eta^{-24}(2\tau)$ is a generator of $A_0(\Gamma_0(2))$ by Lemmata \ref{lem32} and \ref{lem33}, and satisfies the identity \cite{maier}
\begin{align}
	j(\tau) = \dfrac{(f(\tau)+256)^3}{f^2(\tau)}.\label{eq51}
\end{align}
We consider $f(\tau)$ as a modular function on $\Gamma_0(18)$. By Lemma \ref{lem33}, we compute $\mbox{ord}_{1/2} f(\tau)=-9$ and $\mbox{ord}_r f(\tau)=-1$ for $r\in \{\infty, 1/3, 2/3\}$. Since $A_0(\Gamma_0(18)) = \mathbb{C}(f(\tau),u(\tau))$ by Theorem \ref{thm12}, we deduce from Proposition \ref{prop43} that there is a polynomial
\begin{align*}
	F(X,Y)=\sum_{\substack{0\leq i\leq 1\\0\leq j\leq 12}} C_{i,j} X^iY^j\in\mathbb{C}[X,Y]
\end{align*}
with $F(f(\tau), u(\tau))=0$. Substituting the $q$-expansions (\ref{eq41}) of $u(\tau)$ and $f(\tau)$ yields
\begin{align*}
	F(X,Y) = -(1-2Y)^9(1+2Y+4Y^2)+XY(1+Y)^9(1-Y+Y^2),
\end{align*}
so that 
\begin{align}
	f(\tau) = \dfrac{(1-2u)^9(1+2u+4u^2)}{u(1+u)^9(1-u+u^2)},\label{eq52}
\end{align}
where $u:=u(\tau)$. Combining (\ref{eq51}) and (\ref{eq52}) gives 
\begin{align}
	j(\tau) = \dfrac{(1+6u+4u^3)^3P(u)^3}{u(1+u)^9(1-2u)^{18}(1-u+u^2)(1+2u+4u^2)^2},\label{eq53}
\end{align}
where 
\begin{align*}
	P(X):= 1+234X+756X^2+2172X^3+1872X^4+3024X^5+48X^6+3744X^7+64X^9.
\end{align*}
Because $j(\tau)$ is an algebraic integer for $\tau\in K\cap\mathbb{H}$, we see from (\ref{eq53}) that $1/u(\tau)$ is an algebraic integer. Since we know that $1/w(\tau)=1/u(\tau)+1$ from (\ref{eq21}), we conclude that $1/w(\tau)$ is an algebraic integer as well.
\end{proof}

\section{Relations between $u(\tau)$ and $j(d\tau)$}\label{sec6}

We saw in the proof of Theorem \ref{thm52} from Section \ref{sec5} that the integrality of $1/w(\tau)$ at an imaginary quadratic point $\tau\in\mathbb{H}$ follows from expressing the modular $j$-invariant $j(\tau)$ as a rational function of $u:=u(\tau)$, which is given by (\ref{eq53}). Using a similar approach, we derive in this section the formulas for $j(d\tau)$ in terms of $u$ when $d\in \{2,3,6,9,18\}$.

\begin{theorem}\label{thm61}
We have 
\begin{align}
	j(2\tau) = \dfrac{(1+6u^2-2u^3)^3Q(u)^3}{u^2(1-2u)^9(1+u)^{18}(1-u+u^2)^2(1+2u+4u^2)},\label{eq61}
\end{align}
where 
\begin{align*}
	Q(X):=1+234X^2-6X^3+756X^4-936X^5+2172X^6-1512X^7+936X^8-8X^9.
\end{align*}
\end{theorem}

\begin{proof}
We note that $f(\tau) := \eta^{24}(\tau)\eta^{-24}(2\tau)$ satisfies the identity \cite{maier}
\begin{align}
	j(2\tau) = \dfrac{(f(\tau)+16)^3}{f(\tau)}.\label{eq62}
\end{align}
Combining (\ref{eq52}) and (\ref{eq62}) yields (\ref{eq61}).
\end{proof}

\begin{theorem}\label{thm62}
We have 
\begin{align}
	j(3\tau) &= \dfrac{(1+4u^3)^3(1+6u+4u^3)^3(1-6u+36u^2+8u^3-24u^4+16u^6)^3}{u^3(1-2u)^6(1+u)^3(1-u+u^2)^3(1+2u+4u^2)^6},\label{eq63}\\
	j(6\tau) &= \dfrac{(1-2u^3)^3(1+6u-2u^3)^3(1-6u^2-4u^3+36u^4+12u^5+4u^6)^3}{u^6(1-2u)^3(1+u)^6(1-u+u^2)^6(1+2u+4u^2)^3},\label{eq64}\\
	j(9\tau) &= \dfrac{(1+4u^3)^3(1-12u^3+48u^6+64u^9)^3}{u^9(1-2u)^2(1+u)(1-u+u^2)(1+2u+4u^2)^2},\label{eq65}\\
	j(18\tau) &= \dfrac{(1-2u^3)^3(1-6u^3-12u^6-8u^9)^3}{u^{18}(1-2u)(1+u)^2(1-u+u^2)^2(1+2u+4u^2)}.\label{eq66}
\end{align}
\end{theorem}

\begin{proof}
We note that the $\eta$-quotient $g(\tau) := \eta^{12}(\tau)\eta^{-12}(3\tau)$ generates $A_0(\Gamma_0(3))$ by Lemmata \ref{lem32} and \ref{lem33}, and satisfies the identity \cite{maier}
\begin{align}
	j(3\tau) = \dfrac{(g(\tau)+27)(g(\tau)+3)^3}{g(\tau)}.\label{eq67}
\end{align}
We first treat $g(\tau)$ as a modular function on $\Gamma_0(18)$. By Lemma \ref{lem33}, we have $\mbox{ord}_r g(\tau)=-2$ for $r\in \{1/3,2/3,1/9\}$ and $\mbox{ord}_r g(\tau)=-1$ for $r\in \{\infty, 1/6, 5/6\}$. We infer from Theorem \ref{thm12} and Proposition \ref{prop43} that there is a polynomial
\begin{align*}
	G_3(X,Y)=\sum_{\substack{0\leq i\leq 1\\0\leq j\leq 9}} C_{i,j} X^iY^j\in\mathbb{C}[X,Y]
\end{align*}
with $G_3(g(\tau),u)=0$. Using the $q$-expansions (\ref{eq41}) of $u$ and $g(\tau)$, we get 
\begin{align*}
	G_3(X,Y) = (1+Y)^3(1-2Y)^6-XY(1-Y+Y^2)(1+2Y+4Y^2)^2
\end{align*}
and 
\begin{align}
	g(\tau) = \dfrac{(1+u)^3(1-2u)^6}{u(1-u+u^2)(1+2u+4u^2)^2}.\label{eq68}
\end{align}
Combining (\ref{eq67}) and (\ref{eq68}), we arrive at (\ref{eq63}).

We next treat $g(2\tau)\in A_0(\Gamma_0(18))$. We compute, using Lemma \ref{lem33}, $\mbox{ord}_r g(2\tau)=-2$ for $r\in \{\infty, 1/6,5/6\}$ and $\mbox{ord}_r g(2\tau)=-1$ for $r\in \{1/3,2/3,1/9\}$. We deduce from Theorem \ref{thm12} and Proposition \ref{prop43} that there is a polynomial
\begin{align*}
	G_6(X,Y)=\sum_{\substack{0\leq i\leq 1\\0\leq j\leq 9}} C_{i,j} X^iY^j\in\mathbb{C}[X,Y]
\end{align*}
with $G_6(g(2\tau),u)=0$. Using the $q$-expansions (\ref{eq41}) of $u$ and $g(2\tau)$, we have 
\begin{align*}
	G_6(X,Y) = (1+Y)^6(1-2Y)^3-XY^2(1-Y+Y^2)^2(1+2Y+4Y^2)
\end{align*}
and 
\begin{align}
	g(2\tau) = \dfrac{(1+u)^6(1-2u)^3}{u^2(1-u+u^2)^2(1+2u+4u^2)}.\label{eq69}
\end{align}
Substituting (\ref{eq69}) into (\ref{eq67}), we obtain (\ref{eq64}).

We now consider $g(3\tau)\in A_0(\Gamma_0(18))$. We have, by Lemma \ref{lem33}, $\mbox{ord}_{1/9} g(3\tau)=-6$ and $\mbox{ord}_{\infty} g(3\tau)=-1$. We see from Theorem \ref{thm12} and Proposition \ref{prop43} that there is a polynomial
\begin{align*}
	G_9(X,Y)=\sum_{\substack{0\leq i\leq 1\\0\leq j\leq 9}} C_{i,j} X^iY^j\in\mathbb{C}[X,Y]
\end{align*}
with $G_9(g(3\tau),u)=0$. Using the $q$-expansions (\ref{eq41}) of $u$ and $g(3\tau)$, we get
\begin{align*}
	G_9(X,Y) = (1+Y)(1-2Y)^2(1-Y+Y^2)(1+2Y+4Y^2)^2-XY^3
\end{align*}
and 
\begin{align}
	g(3\tau) = \dfrac{(1+u)(1-2u)^2(1-u+u^2)(1+2u+4u^2)^2}{u^3}.\label{eq610}
\end{align}
Substituting (\ref{eq610}) into (\ref{eq67}) gives (\ref{eq65}).

We finally consider $g(6\tau)\in A_0(\Gamma_0(18))$. We have, by Lemma \ref{lem33}, $\mbox{ord}_{1/9} g(6\tau)=-3$ and $\mbox{ord}_{\infty} g(6\tau)=-2$. We see from Theorem \ref{thm12} and Proposition \ref{prop43} that there is a polynomial
\begin{align*}
	G_{18}(X,Y)=\sum_{\substack{0\leq i\leq 1\\0\leq j\leq 9}} C_{i,j} X^iY^j\in\mathbb{C}[X,Y]
\end{align*}
with $G_{18}(g(6\tau),u)=0$. Using the $q$-expansions (\ref{eq41}) of $u$ and $g(6\tau)$, we get
\begin{align*}
	G_{18}(X,Y) = -(1+Y)^2(1-2Y)(1-Y+Y^2)^2(1+2Y+4Y^2)+XY^6
\end{align*}
and 
\begin{align}
	g(6\tau) = \dfrac{(1+u)^2(1-2u)(1-u+u^2)^2(1+2u+4u^2)}{u^6}.\label{eq611}
\end{align}
Combining (\ref{eq611}) and (\ref{eq67}) yields (\ref{eq66}).
\end{proof}

\section*{Acknowledgements}

The authors would like to thank the anonymous referee for insightful comments and for bringing \cite{cooper,hubsy,zagier} to our attention, which simplified our proofs and in general improved the contents of this paper.


\begin{thebibliography}{99}
\bibitem{apery} R. Ap\'{e}ry, Irrationalit\'{e} de $\zeta(2)$ et $\zeta(3)$, {\it Ast\'{e}risque} {\bf 61} (1979) 11--13.

\bibitem{cho} B. Cho, Modular equations for congruence subgroups of genus zero (II), {\it J. Number Theory} {\bf 231} (2022) 48--79.	

\bibitem{chokimkoo} B. Cho, N. M. Kim and J. K. Koo, Affine models of the modular curves $X(p)$ and its application, {\it Ramanujan J.} {\bf 24} (2011) 235--257.

\bibitem{chokoo} B. Cho and J. K. Koo, Construction of class fields over imaginary quadratic fields and applications, {\it Quart. J. Math.} {\bf 61} (2010) 199--216.	

\bibitem{chokoop} B. Cho, J. K. Koo and Y. K. Park, Arithmetic of the Ramanujan--G\"{o}llnitz--Gordon continued fraction, {\it J. Number Theory} {\bf 129} (2009) 922--947.	

\bibitem{coop} S. Cooper, On Ramanujan's function $k(q)=r(q)r(q^2)^2$, {\it Ramanujan J.} {\bf 20} (2009) 311--328.

\bibitem{cooper} S. Cooper, {\it Ramanujan's Theta Functions}, Springer International Publishing AG (2017).

\bibitem{geehons} A. Gee and M. Honsbeek, Singular values of the Rogers-Ramanujan continued fraction, {\it Ramanujan J.} {\bf 11} (2006) 267--284.

\bibitem{hubsy} T. Huber, D. Schultz and D. Ye, Ramanujan-Sato series for $1/\pi$, {\it Acta Arith.} {\bf 207} (2023) 121--160.

\bibitem{ishdai} N. Ishida and N. Ishii, The equations for modular function fields of principal congruence subgroups of prime level, {\it Manuscr. Math.} {\bf 90} (1996) 271--285.

\bibitem{leeparkx} Y. Lee and Y. K. Park, Modular equations of a continued fraction of order six, {\it Open Math.} {\bf 17} (2020) 202--219.

\bibitem{leeparkk} Y. Lee and Y. K. Park, Ramanujan's function $k(\tau)=r(\tau)r(2\tau)^2$ and its modularity, {\it Open Math.} {\bf 18} (2020) 1727--1741.

\bibitem{ligo} G. Ligozat, Courbes modulaires de genus 1, {\it M\'{e}moires de la S. M. F.} {\bf 43} (1975) 5--80.

\bibitem{maier} R. S. Maier, On rationally parametrized modular equations, {\it J. Ramanujan Math. Soc.} {\bf 24} (2009) 1--73, \texttt{https://arxiv.org/abs/math/0611041}.

\bibitem{newman1} M. Newman, Construction and application of a class of modular functions, {\it Proc. London Math. Soc.} {\bf s3--7}  (1957) 334--350.

\bibitem{newman2} M. Newman, Construction and application of a class of modular functions (II), {\it Proc. London Math. Soc.} {\bf s3--9} (1959) 373--387.

\bibitem{ono} K. Ono, {\it The Web of Modularity: Arithmetic of the Coefficients of Modular Forms and $q$-series}, CBMS Regional Conference Series in Mathematics, vol. 102. American Mathematical Society, Rhode Island (2004).

\bibitem{park} Y. K. Park, Analogue of Ramanujan's function $k(\tau)$ for the cubic continued fraction, {\it Int. J. Number Theory} {\bf 19} (2023) 2101--2120.

\bibitem{park2} Y. K. Park, Ramanujan continued fractions of order eighteen, {\it J. Korean Math. Soc.} {\bf 60} (2023) 395--406.

\bibitem{raman} S. Ramanujan, {\it Notebooks (2 volumes)}, Tata Institute of Fundamental Research (1957).

\bibitem{shimura} G. Shimura, {\it Introduction to the Arithmetic Theory of Automorphic Functions}, Princeton University Press (1971).

\bibitem{vasuki} K. R. Vasuki, N. Bhaskar and G. Sharath, On a continued fraction of order six, {\it Ann. Univ. Ferrara} {\bf 56} (2010) 77--89.

\bibitem{ye} D. Ye, Ramanujan's function $k$, revisited, {\it Ramanujan J.} {\bf 56} (2021) 931--952.	

\bibitem{zagier} D. Zagier, Integral solutions of Ap\'{e}ry-like recurrence equations, {\it Groups and Symmetries} (J. Harnad, and P. Winternitz, eds.), CRM Proc. Lecture Notes 47, American Mathematical Society, Rhode Island (2009) 349--366.
\end{thebibliography}
\end{document}